\documentclass[preprint,12pt]{elsarticle}

\usepackage{amssymb}
\usepackage{amsmath}
\input{EMI-macros.sty}

\journal{   }

\begin{document}

\begin{frontmatter}

\title{Analytical solutions for the Extracellular-Membrane-Intracellular model}

\author[label1]{Carlos Ballesteros}
\author[label1]{Alexei Cheviakov}
\author[label2]{Raymond J.~Spiteri}

\affiliation[label1]{organization={Department of Mathematics and Statistics, University of Saskatchewan},
            addressline={106 Wiggins Road}, 
            city={Saskatoon},
            postcode={S7N 5E6}, 
            state={Saskatchewan},
            country={Canada}}
\affiliation[label2]{organization={Department of Computer Science, University of Saskatchewan},
            addressline={110 Science Place}, 
            city={Saskatoon},
            postcode={S7N 5C9}, 
            state={Saskatchewan},
            country={Canada}}

\begin{abstract}
  The Extracellular-Membrane-Intracellular (EMI) model is a novel
  mathematical framework for cardiac electrophysiology
  simulations. The EMI model provides a more detailed description of
  the heart's electrical activity compared to traditional monodomain
  and bidomain models, potentially making it better-suited for
  understanding the electrical dynamics of the heart under
  pathological conditions. In this paper, we derive and verify several
  analytical solutions for the EMI model. Specifically, we obtain a
  family of solutions for a single two-dimensional cell in polar
  coordinates and for a pair of coupled three-dimensional cells in
  spherical coordinates. We also introduce a manufactured solution for
  $N$ three-dimensional cells in Cartesian coordinates. To verify the
  analytical solutions, we conduct numerical experiments using the
  mortar finite element method combined with operator splitting. The
  results demonstrate that the analytical solutions are effective for
  verifying the accuracy of numerical simulations of the EMI model.

\end{abstract}



\begin{keyword}
  Extracellular-Membrane-Intracellular model
  \sep Manufactured solution
  \sep Cardiac electrophysiology
  \sep Operator-splitting method
  \sep Mortar finite element method
\end{keyword}

\end{frontmatter}

\section{Introduction}
\label{sec:Introduction}

Mathematical models of cardiac electrophysiology are highly useful in
understanding the propagation of electrical signals in the heart
tissue.  These models are used to study the effects of various
physiological and pathological conditions, such as cardiac
arrhythmias, on the electrical activity of the heart and to develop
new therapeutic strategies for their treatment. The most popular
cardiac models are the monodomain and bidomain, which work well to
simulate healthy tissue.  However, these models have some limitations
when it comes to simulating an unhealthy heart. To overcome these
limitations, a new model called the
extracellular-membrane-intracellular (EMI) model was
developed~\cite{Tveito2017}. This model considers the tissue as a
collection of individual cells connected to each other, which is more
realistic than the monodomain and bidomain models, where the tissue is
considered as a single homogenized physical region.

Simulating heart tissue at the cellular level using the EMI model 
is computationally intensive due to the complexity and scale 
of the problem. This challenge has driven the development or implementation 
of various numerical methods, especially spatial discretization methods, 
such as the boundary element method~\cite{de2024boundary}, cut finite element
method~\cite{berre2024cut}, mortar finite element method, and
H(div) finite element method~\cite{Tveito2017}.  These implementations
have, in turn, motivated the search for analytical solutions to
validate the accuracy and performance of the numerical methods. To
this end, researchers often opt for the application of manufactured
solutions to simulate a single cell~\cite{Tveito2017,berre2024cut,
ellingsrud2024splitting, fokoue2023numerical}, which provides a useful tool for
convergence analysis. However, these solutions may not fully capture
realistic scenarios.

In this work, we provide a detailed derivation of a family of
analytical solutions for the EMI model that is more general and
realistic than previous analytical solutions for one and two cells. We
also introduce a manufactured solution for $N$ three-dimensional
cells. We numerically implement the EMI model using operator splitting
in conjunction with the mortar finite element method to simulate
simple examples from our solution family for one and two cells, as well
as our manufactured solution, in order to verify the utility of 
the analytical solutions for benchmarking.

The paper is organized as follows. In \Cref{sec:EMI-model}, we provide
a detailed description of the EMI model. In
\Cref{sec:Analytical-solutions}, we derive the family of analytical
solutions. In \Cref{sec:Numerical-solutions}, we describe the
numerical implementation of the EMI model, and in
\Cref{sec:Experiments}, we present the numerical results for some
simple solutions. Finally, in \Cref{sec:Conclusions}, we summarize our
conclusions and outline directions for future work.

\section{The Extracellular-Membrane-Intracellular model}
\label{sec:EMI-model}
The Extracellular-Membrane-Intracellular (EMI) model, also known as 
the cell-by-cell model~\cite{ellingsrud2024splitting, de2024boundary} 
was formally introduced by Tveito et al.~\cite{Tveito2017} in 2017. 
Nonetheless, the model had been the subject of study for many years prior 
to that.~\cite{veneroni2006reaction, franzone2002degenerate, henriquez2017boundary}. 
The model is based on three physical spaces: the \emph{intracellular} space, 
the intracellular-extracellular \emph{membrane} (or cell membrane) space, and 
the \emph{extracellular} space. The tissue is represented as a collection 
of cells spread across the extracellular space that communicate with one 
another via their membranes.

\subsection{Governing equations}
\label{subsec:Governing-equations}
\begin{figure}[H]
    \centering
    \includegraphics[scale=1.1]{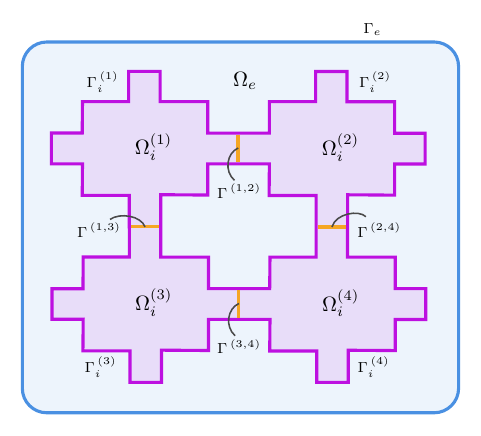}
    \caption{Configuration of a two-dimensional domain with four cells (purple) surrounded by extracellular space (blue) with gap junctions (orange) \cite{Tveito2017}.}
    \label{fig:EMI_Model}
\end{figure}

The complete mathematical description of the EMI model~\cite{Tveito2017, jaeger2021derivation} 
is presented as follows. Let $\Omega_e$ be a bounded domain representing the extracellular space, 
in which $N$ cells are embedded, with its boundary denoted by $\partial\Omegae$. 
For each cell $k$, the intracellular space is represented by the domain $\Omegai[k]$ 
with its boundary $\partial\Omegai[k]$.  The intersection between the boundary of 
cell $k$ and the boundary of the extracellular domain represents the membrane of 
the cell, denoted by $\Gammai[k] := \partial\Omegai[k] \cap \partial\Omegae$, while 
the intersection between the boundaries of two different cells $k$ and $\ell$, 
$\GammaGap[k,\ell] := \partial\Omegai[k] \cap \partial\Omegai[\ell]$, represents a 
gap junction.  The outer boundary of the extracellular space is denoted by $\Gammae$. 
\Cref{fig:EMI_Model} shows a two-dimensional configuration with four cells and 
the associated domains to illustrate the domains of the EMI model.
The intracellular potential $\Omegai[k]$ in the intracellular space of 
cell $k$ is denoted by $\ui[k]$, and the extracellular potential in the 
extracellular space $\Omegae$ is denoted by $\ue$.  The transmembrane 
potential of a cell $k$ is divided into two parts, the transmembrane 
potential across the membrane, $\vk[k]$, which is defined as 
\begin{align}\label{eq:membrane-potential}
    \vk[k] &= \ui[k] - \ue, & \text{on $\Gammai[k]$},
\end{align}
and the transmembrane potential across the gap junctions $\wk[k,\ell]$, 
which for an adjacent cell $l$, is defined as
\begin{align}\label{eq:gap-potential}
    \wk[k,\ell] &= \ui[k] - \ui[\ell] & \text{on $\GammaGap[k,\ell]$}.
\end{align}

In each intracellular space $\Omegai[k]$ and the extracellular domain 
$\Omegae$, the intracellular and extracellular potentials satisfy
\begin{subequations}
  \label{eq:emi-constraints}
  \begin{alignat}{2}
    \nabla\cdot(\sigmai\nabla\ui[k]) &= 0, &&\qquad \text{in $\Omegai[k]$}, \label{eq:emi-intracell-constraint}\\
    \nabla\cdot(\sigmae\nabla\ue) &= 0, && \qquad \text{in $\Omegae$},
    \label{eq:emi-extracell-constraint}
  \end{alignat} 
\end{subequations} 
where the parameters $\sigmai$ and $\sigmae$ represent 
the conductivities of the intracellular and extracellular 
spaces, respectively, where for simplicity we assume $\sigmai$ 
is not cell-dependent. Furthermore, on each membrane $\Gammai[k]$ 
and gap junction $\GammaGap[k,\ell]$, the transmembrane potentials 
$\vk[k]$ and $\wk[k,\ell]$ satisfy \begin{subequations}\label{eq:emi-eqs}
    \begin{align}
        \Cmk[k]\pderiv[]{\vk[k]}{t} &= \Imk[k] - \Iionk[k], &&\text{on $\Gammai[k]$},\label{eq:emi-ap-cellmembrane}\\
        \Imk[k] = (\sigmae\nabla\ue)\cdot\nex &= -(\sigmai\nabla\ui[k])\cdot\nin[k], &&\text{on $\Gammai[k]$},\label{eq:emi-current-cellmembrane}\\
        \Cmk[k,\ell]\pderiv[]{\wk[k,\ell]}{t} &= \Imk[k,\ell] - \Igap[k,\ell], &&\text{on $\GammaGap[k,\ell]$},\label{eq:emi-ap-gapjunction}\\
        \Imk[k,\ell] = (\sigmai\nabla\ui[k])\cdot\nin[k] &= -(\sigmai\nabla\ui[\ell])\cdot\nin[\ell], &&\text{on $\GammaGap[k,\ell]$},\label{eq:emi-current-gapjunction}
    \end{align}
    \end{subequations}
where the parameter $\Cmk[k]$ is the capacitance of the 
membrane $\Gammai[k]$ per unit area, the ion currents $\Imk[k]$ 
and $\Iionk[k]$ are currents across the membrane $\Gammai[k]$, and 
$\nin[k]$ is the outward unit normal vector to $\partial\Omegai[k]$. 
The parameter $\Cmk[k,\ell]$ is the capacitance of the gap junction 
$\GammaGap[k,\ell]$, $\Imk[k,\ell]$ is the ion current across the gap 
junction $\GammaGap[k,\ell]$, and $\Igap[k,\ell]$ is the ion current due 
to the transmembrane potential $\wk[k,\ell]$ at the gap junction. We note 
that at the gap junction $\GammaGap[k,\ell]$, the normal vectors $\nin[k]$ 
from $\partial\Omegai[k]$ and $\nin[\ell]$ from $\partial\Omegai[\ell]$ 
satisfy $\nin[k] = -\nin[\ell]$.

Let $\tn[0],\tend>0$ be such that $\tn[0]<\tend$. The extracellular 
potential $\ue$ is subject to either a Dirichlet boundary condition 
on $\Gammae\times[\tn[0],\tend]$,
\begin{align}\label{eq:emi-dirichletbc}
    \ue &= \uapp,
\end{align}
or a Neumann boundary condition 
\begin{align}\label{eq:emi-neumannbc}
    (\sigmae\nabla \ue)\cdot\nex &= \Iapp.
\end{align}

A cell model is added to the system to couple the transmembrane 
potentials $\vk[k]$ and $\wk[k,\ell]$ with the ion currents $\Iionk[k]$ 
and $\Igap[k,\ell]$.  Cell models fall into two categories: active and 
passive. An active cell model is usually nonlinear and includes multiple 
state variables (see, e.g.,~\cite{ref:tenTusscher2006a,ref:Grandi2009,ref:Gray2016}). 
In contrast, a passive cell model is generally simpler, with the ion 
currents $\Iionk[k]$ and $\Igap[k,\ell]$ assumed to depend linearly on the 
transmembrane potentials $\vk[k]$ and $\wk[k,\ell]$, respectively.  For a 
passive cell model, these ion currents are represented as 
\begin{align}
    \Iionk[k](\vk[k]) &= \frac{1}{\Rmk[k]}\Big(\vk[k]-\vrest\Big),\label{eq:emi-passivecellmodel-membrane}
\end{align}
\begin{align}
    \Igap[k,\ell](\wk[k,\ell]) &= \frac{1}{\Rmk[k,\ell]}\Big(\wk[k,\ell]-\wrest\Big).\label{eq:emi-passivecellmodel-gap}
\end{align}

In case an active cell model is used, the following equation is added 
to the system of equations describing the dynamics of the potentials:
\begin{align}
    \frac{\partial\sk[k]}{\partial t} &= \fk[k](t,\sk[k],\vk[k]), &\text{on $\Gammai[k]$},\label{eq:emi-activecellmodel}
\end{align}
where the vector field $\sk[k]$ describes the cellular state at the 
time $t\in[\tn[0],\tend]$, and the nonlinear vector field $\fk[k]$ 
represents the dynamics of the cell model.
We note that when the passive cell model is implemented, the EMI model 
given by the equations \eqref{eq:membrane-potential}, \eqref{eq:gap-potential}, 
\eqref{eq:emi-constraints}, \eqref{eq:emi-eqs}, \eqref{eq:emi-passivecellmodel-membrane}, 
\eqref{eq:emi-passivecellmodel-gap} subject to \eqref{eq:emi-dirichletbc} or 
\eqref{eq:emi-neumannbc} is a linear system of partial differential equations.
\section{Analytical solutions}
\label{sec:Analytical-solutions}

In this section, analytical solutions are sought for some specific
examples and cases of the EMI model using passive cell models.
Although these solutions may not make physical sense, they provide a
reference point for the validation of numerical methods implemented
to solve real-world problems.

\subsection{Single cell}
\label{subsec:Single-cell}
\begin{figure}[H]
    \centering
    \includegraphics[scale=0.8]{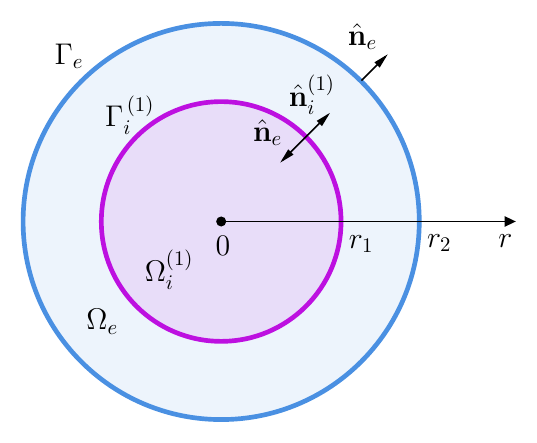}
    \caption{Domain representation for a single circular cell (purple) surrounded by a circular layer as extracellular space (blue).}
    \label{fig:one_spherical_cell}
\end{figure}

Consider a two-dimensional single cell with a circular membrane
surrounded by an outer circular layer representing the extracellular
space, as shown in~\Cref{fig:one_spherical_cell}.  In this setup, no
gap junctions are present. 

Assuming that the transmembrane potential at the cell membrane is set to
$\vk[1] = v_0$ at $t = 0$, we are interested in finding an analytical solution
to the following system of partial differential equations,
\begin{subequations}\label{eq:emipassive}
    \begin{align}
        \nabla\cdot(\sigmai\nabla\ui[1]) &= 0, &&\text{in $\Omegai[1]$},\label{eq:emipassive-intra}\\
        \nabla\cdot(\sigmae\nabla\ue) &= 0, &&\text{in $\Omegae$},\label{eq:emipassive-extra}\\
        \ue &= \uapp, &&\text{on $\Gammae$}.
        \label{eq:emipassive-bc}\\
        \Imk[1] = (\sigmae\nabla\ue)\cdot\nex &= -(\sigmai\nabla\ui[1])\cdot\nin[1], &&\text{on $\Gammai[1]$},\label{eq:emipassive-flux}\\
        \vk[1] &= \ui[1]-\ue, &&\text{on $\Gammai[1]$},\label{eq:emipassive-transpot}\\
        \Cmk[1]\pderiv[]{\vk[1]}{t} &= \Imk[1] - \frac{1}{\Rmk[1]}\Big(\vk[1]-\vrest\Big), &&\text{on $\Gammai[1]$},\label{eq:emipassive-cell}\\
        \vk[1]|_{t=0} &= v_0, &&\text{on $\Gammai[1]$}.\label{eq:emipassive-ic}
    \end{align}
\end{subequations}

Using polar coordinates with the origin at the center of the cell and 
assuming radial symmetry, the potentials $\ui[1]$ and $\ue$ depend 
solely on the radial coordinate $r$.  Solving the equations 
\eqref{eq:emipassive-intra} and \eqref{eq:emipassive-extra}, the
intracellular and extracellular potentials are given by
\begin{subequations}
    \label{eq:sol-single-cell}
\begin{align}
    \ue &= A_{1} \int \frac{1}{\xi \sigmae} \, d\xi + A_{2}, \\
    \ui[1] &= B_{1} \int \frac{1}{\xi \sigmai} \, d\xi + B_{2},
\end{align}
where $A_{1}$, $A_{2}$, $B_{1}$ and $B_{2}$ are the integration constants, 
which are generally arbitrary functions of time.

If we now consider the intracellular and extracellular conductivities 
to be equal, i.e.,
\begin{align}
    \sigmae = \sigmai,
\end{align}
we can find a condition for $A_{1}$ and $B_{1}$ using equation
\eqref{eq:emipassive-flux}:
\begin{align}
    (\sigmae\nabla\ue)\cdot\nex = -(\sigmai\nabla\ui[1])\cdot\nin[1] \Longrightarrow  B_{1} = A_{1} := A,
\end{align}
and, consequently, the transmembrane current and potential take the forms
\begin{align}
    \Imk[1]= (\sigmae\nabla\ue)\cdot\nex = -\frac{A}{r_{1}}, \\
    \vk[1] = \ui[1]-\ue = B_{2}-A_{2},
\end{align}
respectively.

Now that $\Imk[1]$ is known, we can proceed to find the transmembrane 
potential $\vk[1]$ by solving the initial value problem constructed 
with equations \eqref{eq:emipassive-cell} and \eqref{eq:emipassive-ic}. 
The general solution to this system is given by
\begin{align}
    \vk[1] = e^{-\frac{t}{\Cmk[1]\Rmk[1]}}\left(v_0+\int_0^t\frac{e^{\frac{\tau}{\Cmk[1]\Rmk[1]}}}{\Cmk[1]}\Big(\frac{\vrest}{\Rmk[1]}-\frac{A}{r_{1}}\Big)d\tau\right).
\end{align}

Finally, we can find a solution as long as the boundary condition is
satisfied; i.e.,
\begin{align}
    \ue |_{r=r_{2}} = \uapp.
\end{align}
\end{subequations}

In summary, we have a family of solutions where $\sigmae=\sigmai$
(which can be the case for liposomes) and $A$ and $A_{2}$ or $B_{2}$
are free to be chosen.  The solutions presented here can be readily
extended to a three-dimensional space, with a spherical cell and
surrounding shell representing the extracellular space, under the
assumption of radial symmetry.

\subsection{Two coupled semi-spherical cells}
\label{subsec:Two-coupled-cells}

\begin{figure}[htbp]
    \centering  
    \subfigure[Three-dimensional representation of two coupled
    semi-spherical cells.]
    {\resizebox{0.35\textwidth}{!}{\includegraphics[scale=1]{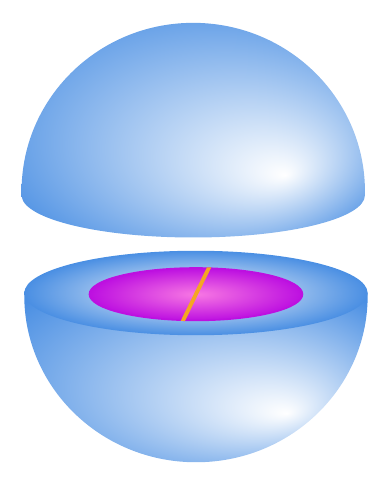}}}
    \hfill
    \subfigure[Two-dimensional cross-section of two semi-spherical cells.]  
    {\resizebox{0.45\textwidth}{!}{\includegraphics[scale=1]{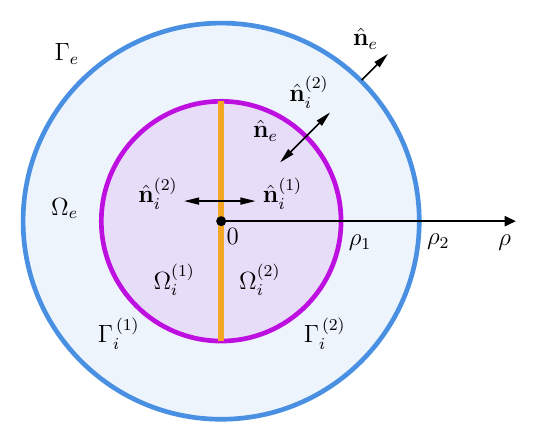}}}
    \caption{Domain representation of two coupled semi-spherical cells (purple) 
    connected by a gap junction (orange) and surrounded by a spherical 
    shell-shaped extracellular space (blue).} 
    \label{fig:two_semispherical_cells}
\end{figure}

Consider two semi-spherical cells connected along their circular bases
and surrounded by a spherical shell of extracellular space as shown
in~\Cref{fig:two_semispherical_cells}.
We assume that the transmembrane potentials at the cell membranes are
initially set to $\vk[k] = v_k$ at $t = 0$ for $k = 1, 2$. Moreover,
the transmembrane potential on the gap junction is initially set to
$\wk[1,2] = w_0$ at $t = 0$.  Our objective now is to find a solution
to the following linear system of partial differential equations:
\begin{subequations}\label{eq:emipassive-twocells}
    \begin{align}
        \nabla\cdot(\sigmai\nabla\ui[k]) &= 0,&&\text{in $\Omegai[k]$},\label{eq:emipassive-intra-twocells}\\
        \nabla\cdot(\sigmae\nabla\ue) &= 0,&&\text{in $\Omegae$},\label{eq:emipassive-extra-twocells}\\
        \ue &= \uapp, &&\text{on $\Gammae$}.\label{eq:emipassive-bc-twocells}\\
        \Imk[k] = (\sigmae\nabla\ue)\cdot\nex &= -(\sigmai\nabla\ui[k])\cdot\nin[k],&&\text{on $\Gammai[k]$},\label{eq:emipassive-flux-twocells}\\
        \vk[k] &= \ui[k]-\ue,&&  \text{on $\Gammai[k]$},\label{eq:emipassive-transpot-twocells}\\
        \Cmk[k]\pderiv[]{\vk[k]}{t} &= \Imk[k] - \frac{1}{\Rmk[k]}\Big(\vk[k]-\vrest \Big),&&\text{on $\Gammai[k]$},\label{eq:emipassive-cell-twocells}\\
        \vk[k]|_{t=0} &= \vk[k]_{0},&&\text{on $\Gammai[k]$},\label{eq:emipassive-ic-twocells}\\
        \Imk[1,2] = (\sigmai\nabla\ui[2])\cdot\nin[2] &= -(\sigmai\nabla\ui[1])\cdot\nin[1],&&\text{on $\GammaGap[1,2]$},\label{eq:emi-current-gapjunction-twocells}\\
        \wk[1,2] &= \ui[1]-\ui[2],&&  \text{on $\GammaGap[1,2]$},\label{eq:emipassive-transpot-w-twocells}\\
        \Cmk[1,2]\pderiv[]{\wk[1,2]}{t} &= \Imk[1,2] - \frac{1}{\Rmk[1,2]}\Big(\wk[1,2]-\wrest\Big),&&\text{on $\GammaGap[1,2]$},\label{eq:emi-gapjunction-twocells}\\
        \wk[1,2]|_{t=0} &= w_{0},&&\text{on $\GammaGap[1,2]$},\label{eq:emipassive-ic-w-twocells}
    \end{align}
\end{subequations}
with $k=1,2$.

We may use the same method as in \Cref{subsec:Single-cell} to obtain
solutions for $\ue$, $\ui[1]$, and $\ui[2]$ by solving equations
\eqref{eq:emipassive-intra-twocells} and \eqref{eq:emipassive-extra-twocells} 
in spherical coordinates and assuming radial symmetry to yield
\begin{subequations}
    \label{eq:sol-two-coupled-cells}
\begin{align}
    \ue &= A_{1} \int \frac{1}{\xi^2 \sigmae} \, d\xi + A_{2}, \\
    \ui[1] &= B_{1} \int \frac{1}{\xi^2 \sigmai} \, d\xi + B_{2},\\
    \ui[2] &= C_{1} \int \frac{1}{\xi^2 \sigmai} \, d\xi + C_{2},
\end{align}
where $A_k$, $B_k$, $C_k$, $k=1,2$, are the constants of integration,
which are generally arbitrary functions of time. Assuming that the
conductivities are equal, i.e.,
\begin{align}
    \sigmai = \sigmae,
\end{align}
equation \eqref{eq:emipassive-flux-twocells}, for $k=1,2$, establishes
the following condition for $A_1$, $B_1$, and $C_1$:
\begin{align}
  B_1 = C_1 = A_1 := A.
\end{align}

Additionally, from equations \eqref{eq:emipassive-transpot-twocells}
and \eqref{eq:emipassive-transpot-w-twocells}, we have
\begin{align}
    \vk[1] &= B_2-A_2,\\
    \vk[2] &= C_2-A_2,
\end{align}
and 
\begin{align}
    \label{eq:w-two-cells-conditions}
    \wk[1,2] = \ui[1]-\ui[2]= B_2-C_2.
\end{align}

The transmembrane current for cell $k$ on the membrane
$\Gammai[k]$ takes the form
\begin{align}
    \Imk[k]= (\sigmai\nabla\ui[k])\cdot\nin[k] = -\frac{A}{(\rho_{1})^2},
\end{align}
while the transmembrane current in the gap junction is
\begin{align}
    \Imk[2,1] =  (\sigmai\nabla\ui[2])\cdot\nin[2] &= -(\sigmai\nabla\ui[1])\cdot\nin[1] = 0.
\end{align}

Now assuming that the cells are characterized by the same model, i.e.,
$\Rmk[1]=\Rmk[2]$ and $\Cmk[1]=\Cmk[2]$, we can define initial-value
problems using equations \eqref{eq:emipassive-cell-twocells} and
\eqref{eq:emipassive-ic-twocells} to solve for potentials
$\vk[k]$, $k=1,2$, resulting in
\begin{align}
    \vk[k] = e^{-\frac{t}{\Cmk[1]\Rmk[1]}}\left(\vk[k]_0+\int_0^t\frac{e^{\frac{\tau}{\Cmk[1]\Rmk[1]}}}{\Cmk[1]}\Big(\frac{\vrest}{\Rmk[1]}-\frac{A}{(\rho_{1})^2}\Big)d\tau\right).
\end{align}

From \cref{eq:w-two-cells-conditions}, $\wk[1,2]$ takes the form
\begin{align}
    \wk[2,1]=(\vk[1]_0-\vk[2]_0) e^{-\frac{t}{\Cmk[1]\Rmk[1]}}.
\end{align}

We note that the expression above satisfies the initial-value problem
created with equations \eqref{eq:emi-gapjunction-twocells} and
\eqref{eq:emipassive-ic-w-twocells} in two cases. The first case is
when $\vk[2]_0=\vk[1]_0$ and $\wrest=w_0=0$; the second case
is when $\vk[2]_0\neq \vk[1]_0$, $w_0=0$, and
$\Cmk[1,2]\Rmk[1,2]=\Cmk[1]\Rmk[1]$. Finally, we are able to determine
a solution subject to satisfaction of the boundary condition 
\begin{align}
    \ue |_{\rho=\rho_{2}} = \uapp.
\end{align}
\end{subequations}

In summary, we have derived a family of solutions under the following
conditions: the conductivities $\sigmai$ and $\sigmae$ are equal, both
cells are identical, the resting potential of the gap junction model
is zero (which is realistic according to~\cite{Tveito2017}), and the
ratio of the membrane capacitance to the gap junction capacitance is
equal to the ratio of the cell membrane resistance to the gap junction
resistance. Additionally, $A$ and any one of $A_2$, $B_2$,
or $C_2$ are free to be chosen.

\subsection{Manufactured solutions}
\label{subsec:Manufactured-solutions}

The method of manufactured solutions (MMS) is a commonly used 
approach to construct analytical solutions for the differential 
equations that underlie a simulation code. The key idea behind MMS 
is to select a smooth and sufficiently differentiable function that 
satisfies the differential equations and boundary conditions of 
the problem. The differential equation is then modified by adding 
forcing (or \emph{source}) terms to ensure that the chosen function 
becomes a solution of the modified equation. This technique has been 
used to validate numerical methods for solving the EMI model of a 
single cell; see, e.g.,~\cite{Tveito2017, ellingsrud2024splitting, berre2024cut}.  
In this subsection, we present a manufactured solution for $N$ 
three-dimensional cells.

\begin{figure}[htbp]
    \centering
    \subfigure[$xy$-plane representation of the single cell.]
    {\resizebox{0.45\textwidth}{!}{\includegraphics[scale=2]{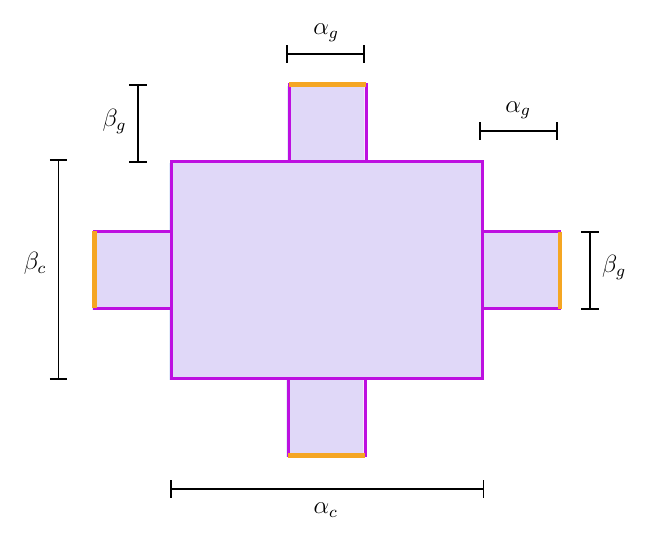}}}
    \hfill
    \subfigure[$yz$-plane representation of the single cell.]  
    {\resizebox{0.4\textwidth}{!}{\includegraphics[scale=2]{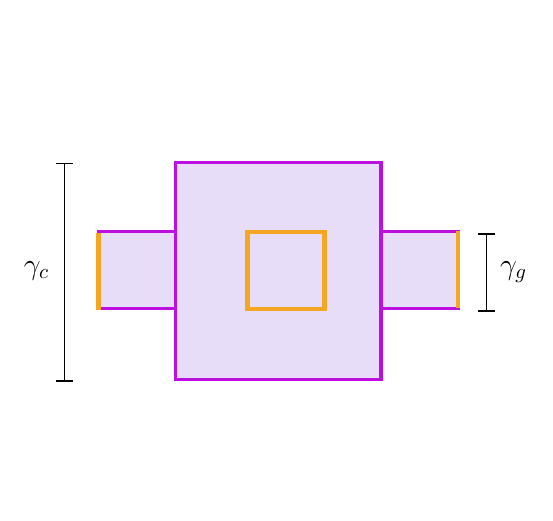}}}
    \caption{Single three-dimensional cardiac cell and its dimensions:
      cell membranes (purple) and gap junctions (orange).}
    \label{fig:manufactured-solutions}
\end{figure}

Consider a three-dimensional cuboidal extracelullar domain with $N$
identical three-dimensional cells embedded in it arranged in a
sheet-like structure. The cells are connected with each other through
gap junctions; see \Cref{fig:manufactured-solutions} for a depiction
and dimensions of a single cell. 
To construct a manufactured solution, we require $\ue$ and $\ui[k]$ to
satisfy Poisson's equation, that is,
\begin{subequations}\label{eq:sol-manufactured}
\begin{align}
     -\nabla\cdot(\sigmae\nabla\ue) &= \fe \label{eq:poisson_equations_ue},\\
     -\nabla\cdot(\sigmai\nabla\ui[k])&=\fik[k] \label{eq:poisson_equations_ui},
\end{align}
with $k=1,2,\dots,N$. 

Using the method of separation of variables, we propose the 
following form for the extracellular
potential $\ue$ and the intracellular potentials $\ui[k]$:

\begin{align}
    \ue &= \Te(t)X(x,y,z), && (x,y,z,t)\in\Omegae\times[\tn[0],\tend],\\
    \ui[k] &= \Ti[k](t)X(x,y,z), && (x,y,z,t)\in\Omegai[k]\times[\tn[0],\tend],
\end{align}
where $\Te=\Te(t)$, $\Ti[k]=\Ti[k](t)$, $k=1,2,\dots,N$, and
$X=X(x,y,z)$ are functions to be determined.

The continuity condition for the normal fluxes 
in \eqref{eq:emi-current-cellmembrane} and \eqref{eq:emi-current-gapjunction} becomes
\begin{align*}
    (\sigmai\Ti[k]-\sigmae\Te)\nabla X\cdot\nin[k] &= 0, && \text{on $\Gammai[k]$},
\end{align*}
and 
\begin{align*}
    \sigmai(\Ti[k]-\Ti[\ell])\nabla X\cdot\nin[k]& = 0, && \text{on $\GammaGap[k,\ell]$}.
\end{align*}

To ensure the continuity of normal fluxes across cell membranes and gap 
junctions, we impose the continuity of the normal flux of $X$ at these domains, that is
\begin{align*}
    \nabla X\cdot \hatn &= 0, && \text{on $\Gammai[k]\cup\GammaGap[k,\ell]$},
\end{align*}
where $\hatn$ is the unit normal vector on
$\Gammai[k]\cup\GammaGap[k,\ell]$, which is defined as $\hatn = \nin[k]$
on $\Gammai[k]$ and $\hatn = \nin[k]$ on $\GammaGap[k,\ell]$. Now
if we set the origin of the system at the center of any cell and
consider the dimensions $\alpha_g$, $\beta_g$,$\gamma_g$, $\alpha_c$,
$\beta_c$, and $\gamma_c$ described
in~\Cref{fig:manufactured-solutions}, we can choose the following form
for the function $X$:
\begin{align}
    X &= \cos\left(\frac{2\pi x}{\alpha_g}\right)\cos\left(\frac{2\pi y}{\beta _g}\right)
        \cos\left(\frac{2\pi z}{\gamma_g}\right),\label{eq:space-func}
\end{align}
where $\alpha_g$, $\beta_g$, and $\gamma_g$ are divisors of $\alpha_c$, 
$\beta_c$, and $\gamma_c$, respectively.

We see now that this choice implies that $\Imk[k] = 0$ on
$\Gammai[k]$ and that $\Imk[k,\ell] = 0$ on $\GammaGap[k,\ell]$. Therefore,
using the passive cell models \eqref{eq:emi-passivecellmodel-membrane}
and \eqref{eq:emi-passivecellmodel-gap}, the partial differential equations 
for the membrane potentials $\vk[k]$ and the gap junctions $\wk[k,\ell]$ in 
\eqref{eq:emi-eqs} take the form

\begin{align}
    \Cmk[k]\pderiv[]{\vk[k]}{t} + \frac{1}{\Rmk[k]}\Big(\vk[k]-\vrest\Big) &= \gk[k], && \text{on $\Gammai[k]$},\\
    \Cmk[k,\ell]\pderiv[]{\wk[k,\ell]}{t} + \frac{1}{\Rmk[k,\ell]}\Big(\wk[k,\ell]-\wrest\Big) &= \gk[k,\ell], && \text{on $\Gammai[k,\ell]$},
\end{align}
where $\gk[k]$ and $\gk[k,\ell]$ are forcing functions to be determined.

We can now choose the form of the time-dependent functions $\Ti[k]$ and $\Te$.
Using a similar approach as for $X$, we select the following simple solutions
\begin{align}
    \Ti[k] &= \Ak[k]e^{\left(-\frac{t}{\Cmk[k]\Rmk[k]}\right)}+\Te,\quad \Te = B,\label{eq:time-funcs}
\end{align}
where $\Ak[k]$, and $B$ are constants that represent the initial values 
of the functions $\Ti[k]$, and $\Te$. 

Using \eqref{eq:membrane-potential} and \eqref{eq:gap-potential}, we can define
$\gk[k]$ and $\gk[k,\ell]$ as
\begin{align}
    \gk[k] &= -\frac{\vrest}{\Rmk[k]},\\
    \gk[k,\ell] &= \left(\frac{1}{\Rmk[k,\ell]}-\frac{\Cmk[k,\ell]}{\Cmk[k]\Rmk[k]}\right)\ui[k]-\left(\frac{1}{\Rmk[k,\ell]}-\frac{\Cmk[k,\ell]}{\Cmk[\ell]\Rmk[\ell]}\right)\ui[\ell] - \frac{\wrest}{\Rmk[k,\ell]},
\end{align}
respectively.

Finally, using \eqref{eq:poisson_equations_ue} and \eqref{eq:poisson_equations_ui} 
the forcing functions $\fe$ and $\fik[k]$ can be taken as
\begin{align}
    \fe &= 
        \sigmae\left(\left(\frac{2\pi}{\alpha_g}\right)^2+\left(\frac{2\pi}{\beta_g}\right)^2+\left(\frac{2\pi}{\gamma_g}\right)^2\right)\ue,\\
    \fik[k] &= 
        \sigmai\left(\left(\frac{2\pi}{\alpha_g}\right)^2+\left(\frac{2\pi}{\beta_g}\right)^2+\left(\frac{2\pi}{\gamma_g}\right)^2\right)\ui[k],
\end{align}
and the boundary condition on $\Gammae\times[\tn[0],\tend]$ can be
given by
\begin{align}
    \uapp &=\Te X,
\end{align}
for Dirichlet boundary condition and
\begin{align}
    \Iapp &= \Te\sigmae\nabla X\cdot\nex,
\end{align}
for Neumann boundary condition.
\end{subequations}

\section{Numerical solutions}
\label{sec:Numerical-solutions}
This section describes a numerical method to simulate the presented
exact solutions to the EMI model.  We employ the mortar finite element
method to discretize the spatial domain and operator splitting to
discretize the time evolution.

\subsection{Spatial discretization using the mortar finite element method}
\label{subsec:Spatial-discretization-using-mortar-finite-element-method}

For simplicity, we present the spatial discretization 
for two cells using \Cref{eq:emipassive-twocells} 
subject to Dirichlet boundary conditions. This discretization can
easily be extended to $N$ cells or reduced to a single cell. Our
presentation follows a similar approach to that in~\cite{Tveito2017}.

Multiplying by an appropriate test function $\phi_i$ and integrating by parts, 
equations \eqref{eq:emipassive-intra-twocells} and
\eqref{eq:emipassive-extra-twocells} for $k=1,2$, can be
rewritten in the weak form as
\begin{subequations}\label{eq:emi-varproblem}
\begin{align}
    \int_{\Omegai[1]} \sigmai \nabla \ui[1] \cdot \nabla \phi_i^{(1)} \, dx + \int_{\Gammai[1]} \Imk[1] \phi_i^{(1)} \, ds + \int_{\GammaGap[2,1]} \Imk[1,2] \phi_i^{(1)} \, ds &= 0 && \forall \phi_i^{(1)} \in V_i^{(1)}, \\
    \int_{\Omegai[2]} \sigmai \nabla \ui[2] \cdot \nabla \phi_i^{(2)} \, dx + \int_{\Gammai[2]} \Imk[2] \phi_i^{(2)} \, ds - \int_{\GammaGap[2,1]} \Imk[1,2] \phi_i^{(2)} \, ds &= 0 && \forall \phi_i^{(2)} \in V_i^{(2)}, \\
    \int_{\Omegae}\sigmae \nabla \ue \cdot \nabla \phi_e \, dx - \int_{\Gammai[1]} \Imk[1] \phi_e \, ds - \int_{\Gammai[2]} \Imk[2] \phi_e \, ds &= 0 && \forall \phi_e \in V_e,
\end{align}
where $V_{e}$ and $V_{i}^{(k)}$ are the function spaces on $\Omegae$
and $\Omegai[k]$, respectively.

Let $G$ be a function space defined over the cell membranes,
specifically
$\Gamma=\Gammai[1]\cap\Gammai[2]\cap\GammaGap[1,2]$. Multiplying
equation \cref{eq:emipassive-cell-twocells} by the test
function $\psi^{(k)}$, associated with the space $\Gammai[k]$ for
$k=1,2$, and equation \eqref{eq:emi-gapjunction-twocells} by
$\psi^{(1, 2)}$, associated with the space $\GammaGap[1,2]$, and
integrating over their respective domains, the equations take the form
\begin{align}
    \int_{\Gammai[1]}\Cmk[1] \pderiv[]{\vk[1]}{t} \psi^{(1)}\, ds &= \int_{\Gammai[1]}\Imk[1]\psi^{(1)} - \frac{1}{\Rmk[1]}\big(\vk[1]-\vrest \big)\psi^{(1)} \, ds && \forall \psi \in G,\\
    \int_{\Gammai[2]}\Cmk[2] \pderiv[]{\vk[2]}{t} \psi^{(2)}\, ds &= \int_{\Gammai[2]}\Imk[2]\psi^{(2)}- \frac{1}{\Rmk[2]}\big(\vk[2]-\vrest \big)\psi^{(2)} \, ds  && \forall \psi \in G,\\
    \int_{\GammaGap[1,2]}\Cmk[1,2] \pderiv[]{\wk[1,2]}{t} \psi^{(1,2)}\, ds &= \int_{\GammaGap[1,2]}\Imk[1,2]\psi^{(1,2)} - \frac{1}{\Rmk[1,2]}\big(\wk[1,2]-\wrest\big)\psi^{(1,2)}\, ds && \forall \psi \in G.
\end{align}

Using the same function space defined above, we require our
variational problem to satisfy the equations
\eqref{eq:emipassive-transpot-twocells} and
\eqref{eq:emipassive-transpot-w-twocells}. This is achieved
by multiplying each equation by the associated test function
$\gamma^{j}$, where $j$ corresponds to the respective space as
previously defined, and integrating over the appropriate domain. That
is,
\begin{align}
    \int_{\Gammai[1]} (\vk[1] - \ui[1] + \ue) \gamma^{(1)} \, ds &= 0 && \forall \gamma \in G,\\
    \int_{\Gammai[2]} (\vk[2]-\ui[2] + \ue)  \gamma^{(2)} \, ds &= 0 && \forall \gamma \in G,\\
    \int_{\GammaGap[1,2]} (\wk[1,2]-\ui[1] + \ui[2]) \gamma^{(1,2)} \, ds &= 0 && \forall \gamma \in G.
\end{align}

\end{subequations}

To describe the finite element discretization of the well-posed problem
formulated using equations \eqref{eq:emi-varproblem}, we introduce
$T_{e,h}$ and $T_{i,h}^k$ as the meshes defined over the domains
$\Omegai[k]$ and $\Omegae$, respectively. The discrete finite element
subspaces $V_{e,h}$, $V_{i,h}^k$, and $G_h$, with $h$ representing the
mesh size parameter, are constructed using continuous piecewise linear
Lagrange elements. Each subspace is defined on its respective mesh,
and the functions in these spaces are linear polynomials on each
element of the mesh. That is,
\begin{align*}
V_{e,h} &= \{v \in C(T_{e,h}): v|_K = P^1(K) \quad \forall K \in T_{e,h}\}, \\
V_{i,h}^k &= \{v \in C(T_{i,h}^k): v|_K = P^1(K) \quad \forall K \in T_{i,h}^k\}, \\
G_h &= \{g \in C(T_h)\phantom{^k}: g|_K = P^1(K) \quad \forall K \in T_h\}.
\end{align*}

\subsection{Time evolution via operator splitting}
\label{subsec:Time-evolution-via-operator-splitting-methods}

The main idea behind operator splitting is to improve efficiency by
splitting a differential equation into simpler parts, such as its
linear and nonlinear parts or its stiff and non-stiff parts. This
allows each part to be solved with methods that suit its specific
characteristics. For example, consider an initial-value problem (IVP)
\begin{subequations}\label{eqn:op-subproblems}
  \begin{gather}
    \frac{\mathrm{d} \yy}{\mathrm{d} t} = \ff(t,\yy):= \ff^{[1]}(t,\yy)+\ff^{[2]}(t,\yy),  \quad \yy(t_{0})=\yy_{0}, \label{eqn:op-subproblems-a}\\
    \frac{\mathrm{d} \yy^{[1]}}{\mathrm{d} t} =
    \ff^{[1]}(t,\yy^{[1]}), \quad \quad \frac{\mathrm{d}
      \yy^{[2]}}{\mathrm{d} t} =
    \ff^{[2]}(t,\yy^{[2]}), \label{eqn:op-subproblems-b}
  \end{gather}
\end{subequations}
where the right-hand side function $\ff(t,\yy)$ of the original
ordinary differential equation (ODE)~\cref{eqn:op-subproblems-a} is
additively split into two parts (or \textit{operators}),
$\fopk[1](t,\yy)$ and $\fopk[2](t,\yy)$, which are then used to form
two ODEs~\cref{eqn:op-subproblems-b}.

To now illustrate the workings of operator-splitting methods, consider
the simplest approach: the Lie--Trotter (or Godunov)
method~\cite{trotter1958approximation, godunov1959difference}.
Applied to the IVP \cref{eqn:op-subproblems-a}, this method begins by
solving the subproblem defined by the operator $\ff^{[1]}$ over
a time step $\Delta t$, using the initial condition from the original
problem.  The solution of this first subproblem, $\yy^{[1]}$,
then serves as the initial condition for the second subproblem defined
by $\ff^{[2]}$, which is also solved over $\Delta t$. The
solution of the second subproblem provides the approximate solution
for the original problem at the end of the time step.

Operator-splitting methods applied to cardiac simulations typically 
separate the problem into two parts: the cell model, which may be 
stiff (e.g.,~\cite{ref:tenTusscher2006a, ref:Grandi2009}), and the 
spatial propagation (or remainder of the system). However, is it 
possible to split the problem more, for example, by separating the
gating equations as well (see, e.g.,~\cite{green2019gating}). 
In this work, we adopt a two-splitting of the EMI model, separating 
the cell models from the propagation component, following 
Tveito et al.~\cite{Tveito2017} and Dominguez et 
al.~\cite{dominguez2021simulation}. Yet, a three-splitting, 
where the membrane and gap junctions are treated separately, naturally 
arises in this context and may offer advantages, particularly when 
evaluating high-order methods~\cite{spiteri2023beyond}.

For applying the operator-splitting method to equations
\eqref{eq:emi-varproblem}, we can consider the problem as
\begin{align*}
  \frac{\partial \Ae \yye}{\partial t} = \ffe^{[1]}(\yye)+\ffe^{[2]}(\yye),  \quad \yye(t_{0})=\yyeo,\\
\end{align*}
where $\yye = (\ue, \ui[1], \ui[2], \vk[1], \vk[2], \wk[1,2])^T$,
$\Ae$ is a matrix such that
$\Ae \yye = (0, 0, 0, \vk[1], \vk[2], \wk[1,2], 0, 0, 0)^T$, and
$\yyeo$ represents the initial data for variables in
$\yye$. Therefore, $\ffe^{[1]}$ and $\ffe^{[2]}$ take the form
\begin{align*}
  \ffe^{[1]}(\yye) &= \begin{pmatrix}
    0\\
    0\\
    0\\[1em]
    - \frac{1}{\Cmk[1]\Rmk[1]}(\vk[1]-\vrest)\\[1em]
    - \frac{1}{\Cmk[2]\Rmk[2]}(\vk[2]-\vrest)\\[1em]
    - \frac{1}{\Cmk[1,2]\Rmk[1,2]}(\wk[1, 2]-\wrest)\\[1em]
    0\\
    0\\
    0\\
  \end{pmatrix},\\
\end{align*}

\begin{align*}
  \ffe^{[2]}(\yye) &= \begin{pmatrix}
    \int_{\Omegai[1]} \sigmai \nabla \ui[1] \cdot \nabla \phi_i^{(1)} \, dx + \int_{\Gammai[1]} \Imk[1] \phi_i^{(1)} \, ds + \int_{\GammaGap[2,1]} \Imk[1,2] \phi_i^{(1)} \, ds \\[1em]
    \int_{\Omegai[2]} \sigmai \nabla \ui[2] \cdot \nabla \phi_i^{(2)} \, dx + \int_{\Gammai[2]} \Imk[2] \phi_i^{(2)} \, ds - \int_{\GammaGap[2,1]} \Imk[1,2] \phi_i^{(2)} \, ds \\[1em]
    \int_{\Omegae}\sigmae \nabla \ue \cdot \nabla \phi_e \, dx - \int_{\Gammai[1]} \Imk[1] \phi_e \, ds - \int_{\Gammai[2]} \Imk[2] \phi_e \, ds \\[1em]
    \int_{\Gammai[1]}\frac{\Imk[1]}{\Cmk[1]} \psi^{(1)} \, ds \\[1em]
    \int_{\Gammai[2]}\frac{\Imk[2]}{\Cmk[2]} \psi^{(2)} \, ds \\[1em]
    \int_{\GammaGap[1,2]}\frac{\Imk[1,2]}{\Cmk[1,2]}\psi^{(1,2)}\, ds \\[1em]
    \int_{\Gammai[1]} (\vk[1] - \ui[1] + \ue) \gamma^{(1)} \, ds \\[1em]
    \int_{\Gammai[2]} (\vk[2]-\ui[2] + \ue)  \gamma^{(2)} \, ds \\[1em]
    \int_{\GammaGap[1,2]} (\wk[1,2]-\ui[1] + \ui[2]) \gamma^{(1,2)} \, ds \\[1em]
  \end{pmatrix}.\\
\end{align*}

This methodology was implemented using the \texttt{Firedrake}
package~\cite{FiredrakeUserManual} for spatial discretization, and the
\texttt{pythOS} and \texttt{irksome}
libraries~\cite{guenter2024pythos, farrell2021irksome} to handle time
evolution via operator splitting.

\section{Examples of numerical solutions}
\label{sec:Experiments}
In this section, we present examples of numerical computations that were
performed in order to verify the agreement with and stability of the 
exact solutions to the EMI model as described in 
\Cref{sec:Analytical-solutions} using the numerical approach described 
in \Cref{sec:Numerical-solutions}.

\subsection{Accuracy}
To assess the accuracy of the exact solutions described in
\cref{subsec:Time-evolution-via-operator-splitting-methods}, we
evaluate the $L^2$-norm of the error for each variable $\vk[k]$,
$\ui[k]$, $\wk[k,\ell]$, and $\ue$ in experiments. The $L^2$-norm of the
error of a quantity $y$ is defined as
\begin{align*}
    \|\tilde{y}-y\|_{2} = \left( \int_\Omega |\tilde{y}-y|^2 \, dx \right)^{\frac{1}{2}},
\end{align*}
where $\tilde{y}$ and $y$ denote the numerical and analytical
solutions, respectively, and $\Omega$ represents the domain of $y$.

\subsection{Experiments}
To simulate the examples presented in
\cref{subsec:Single-cell,subsec:Two-coupled-cells}, it is necessary to
exclude the origin to avoid the singularity inherent in the solutions
and ensure convergence of the numerical approximation. Accordingly, we
exclude a circular or hemisphere region centered at the origin within
the respective domains and impose Dirichlet boundary conditions on the
resulting boundary for the corresponding variables. The meshes for
these domains are generated using
\texttt{Gmsh}~\cite{geuzaine2009gmsh}.

For simplicity, we set the EMI model parameters as
$\Cmk[k] = \Cmk[k,\ell] = 1$ for all simulations. For Experiments 1 and
2, the boundary has a radius of $6$, while the intracellular space has
a radius of $5$. Additionally, the origin is excluded by removing a
circular or spherical region of radius $3$. The conductivities are set
to $\sigmai = \sigmae = 1$.

\subsubsection{Experiment 1: Sinusoidal pulse}
\label{subsubsec:Experiment-1}
Choosing the cell model with parameters $\Rmk[1]=1$ and $\vrest=5$, and
setting $A=10 \sin(t)$, $A_2=5$ and $v_0 = 20$ in equations 
\eqref{eq:sol-single-cell}, the exact solution for a single cell takes the form of a
sinusoidal pulse:
\begin{align*}
        \ue &= 10 \sin(t) \ln(r) + 5,\\
        \ui[1] &= \sin(t) \big( 10\ln(r) -1 \big)+\cos(t)+14 e^{-t}+10,\\
        \vk[1] &= 14 e^{-t} + \cos(t) - \sin(t) + 5,
\end{align*}
    with extracellular boundary condition
\begin{align*}
    \uapp = 10 \ln(6)\sin(t) + 5.
\end{align*}

\setlength{\tabcolsep}{10pt}
\renewcommand{\arraystretch}{1.2}
\begin{table}[H]
  \caption{Errors for Experiment 1. Simulation time spans are from
    $\tn[0]=0.25$ to $\tend=7$, $\nf$ denotes the number of time
    steps, and $\cl$ represents the characteristic length of the
    mesh.}
  \label{tab:experiment-1}
\begin{center}
\begin{adjustbox}{width=0.6 \textwidth, center}
  \begin{tabular}{c c c c c} \hline
   \rule{0pt}{18pt}
    $\cl$ & $\nf$ & $\left\| \tilde{u}_{e} - u_{e}\right\|_{2}$ & 
        $\left\| \tilde{u}^{(1)}_{i} - \ui[1]\right\|_{2}$ & 
            $\left\| \tilde{v}^{(1)}-\vk[1] \right\|_{2}$ \\[7pt] \hline
    0.40  &  7   & \num{9.035e-1} & \num{3.575e0} & \num{5.929e0}  \\
    0.28  &  14  & \num{4.164e-1} & \num{1.651e0} & \num{2.736e0}  \\
    0.20  &  28  & \num{2.016e-1} & \num{8.006e-1} & \num{1.326e0}  \\
    0.14  &  56  & \num{9.937e-2} & \num{3.947e-1} & \num{6.536e-1}  \\
    0.10  &  112 & \num{4.932e-2} & \num{1.960e-1} & \num{3.246e-1}  \\ \hline
  \end{tabular}
\end{adjustbox}
\end{center}
\end{table}

\begin{figure}[H]
  \centering
  \includegraphics[width=1\textwidth]{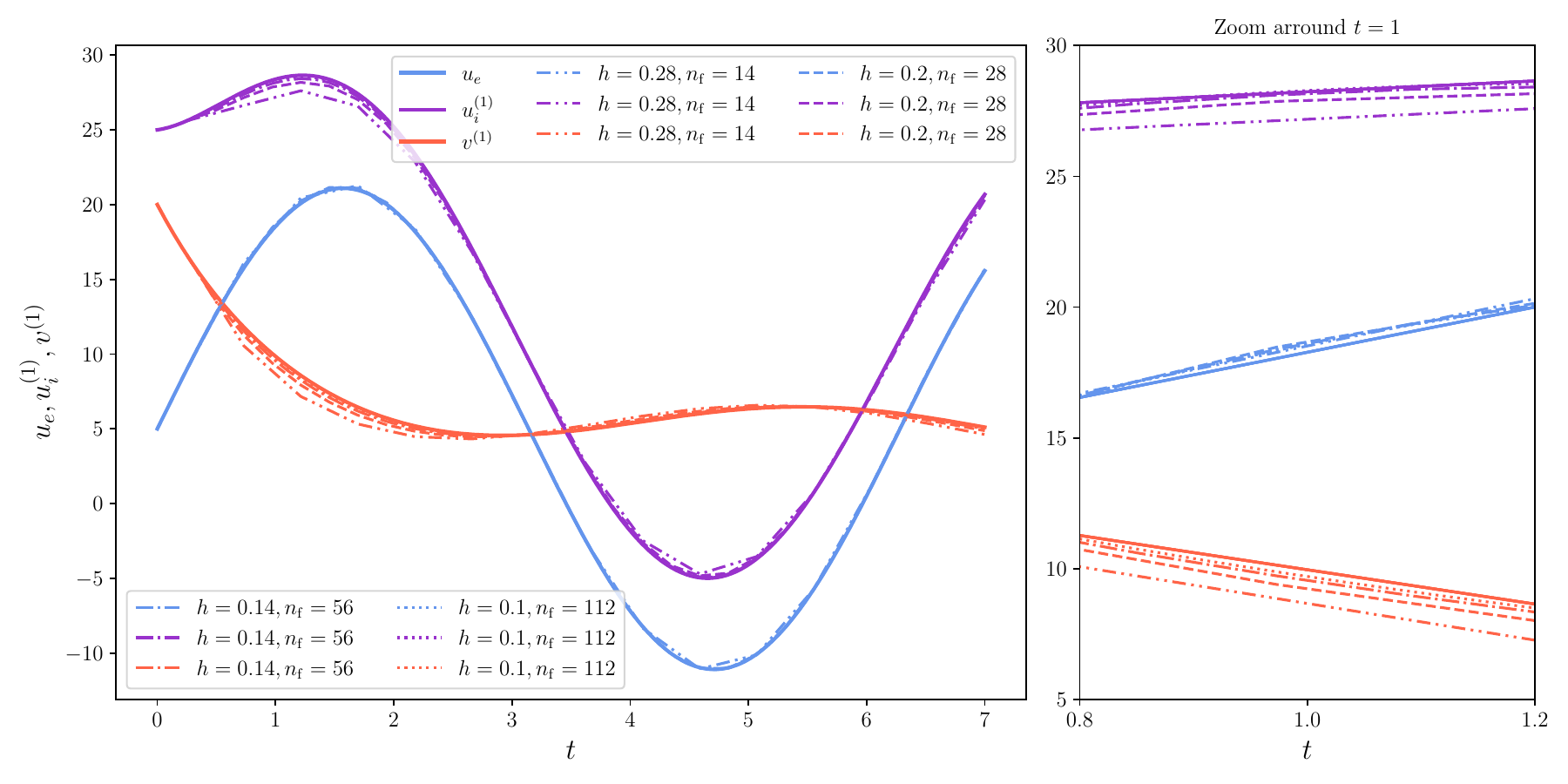}
  \caption{Numerical and exact solutions for $\ue$, $\ui[1]$, and
    $\vk[1]$ at $r=5$, in Experiment 1.}
    \label{fig:plot1_exp1}
\end{figure}

\Cref{tab:experiment-1} demonstrates the correct implementation 
of the numerical method for the single-cell model. The table 
shows that the $L^2$-norm of the error for the variables 
$\vk[1]$, $\ui[1]$, and $\ue$ decreases as the mesh size 
and time step decrease. In addition, \cref{fig:plot1_exp1} illustrates 
that the numerical solutions converge to the exact solutions as 
the mesh size and time step are refined.

\subsubsection{Experiment 2: Damped cosine pulse}
\label{subsubsec:Experiment-2}
Using the same cell model as in \cref{subsubsec:Experiment-1} and
setting $A=-50 e^{-\frac{t}{10}}\cos(t)$, $A_2=0$, $v_{0}^{(1)} = 10$
and $v_{0}^{(2)} = 30$ in equations \eqref{eq:sol-two-coupled-cells},
the exact solution for two-coupled spherical cells takes the form of a
damped cosine pulse:
\begin{align*}
      \ue &= \frac{50 e^{-\frac{t}{10}} \cos(t)}{\rho},\\
      \ui[k] &= 5 + \frac{e^{-t}(144k + 581)}{181} + \frac{2 e^{-\frac{t}{10}} \big( (905 + 18\rho) \cos(t) + 20 \rho \sin(t) \big)}{181 \rho},\\
      \vk[k] &= 5 + \frac{e^{-t}(144k + 581)}{181} + \frac{20 e^{-\frac{t}{10}} \big( 9 \cos(t) + 10 \sin(t) \big)}{181},\\
      \wk[1,2] &= -20 e^{-t},
\end{align*}
    with $k=1,2.$ The extracellular boundary condition is given by
\begin{align*}
    \uapp = \frac{50 e^{-\frac{t}{10}} \cos(t)}{6}.
\end{align*}

\setlength{\tabcolsep}{10pt}
\renewcommand{\arraystretch}{1.2}
\begin{table}[H]
  \caption{Errors for Experiment 2. Simulation time spans are from
    $\tn[0]=0.25$ to $\tend=7$, $\nf$ denotes the number of time
    steps, and $\cl$ represents the characteristic length of the
    mesh. The calculation domain for variables $\ui[1]$, $\ui[2]$, and
    $\wk[1,2]$ is $4\leq \rho \leq5$.  }
  \label{tab:experiment-2}
\begin{center}
\begin{adjustbox}{width=1 \textwidth, center}
  \begin{tabular}{c c c c c c c c } \hline
   \rule{0pt}{18pt}
   $\cl$ & $\nf$ & $\left\| \uet - \ue \right\|_{2}$ & 
   $\left\| \uit[1] - \ui[1] \right\|_{2}$ & 
   $\left\| \uit[2] - \ui[2] \right\|_{2}$ &
   $\left\| \vkt[1] - \vk[1] \right\|_{2}$ &
   $\left\| \vkt[2] - \vk[2] \right\|_{2}$ &
   $\left\| \wkt[1,2] - \wk[1,2] \right\|_{2}$ \\[7pt] \hline
      0.40  &  10   & \num{5.107e-1} & \num{3.311e0} & \num{2.987e0} & \num{4.476e0} & \num{4.361e0} & \num{3.921e-1}  \\
      0.28  &  20  & \num{1.985e-1} & \num{1.872e0} & \num{1.457e0} & \num{2.237e0} & \num{2.061e0} & \num{3.224e-1} \\
      0.20  &  40  & \num{8.778e-2} & \num{1.129e0} & \num{7.751e-1} & \num{1.211e0} & \num{1.055e0} & \num{2.499e-1}  \\
      0.14  &  80  & \num{4.840e-2} & \num{7.161e-1} & \num{4.362e-1} & \num{6.954e-1} & \num{5.672e-1} & \num{1.887e-1} \\
      0.10  &  160  & \num{3.404e-2} & \num{4.701e-1} & \num{2.573e-1} & \num{4.182e-1} & \num{3.174e-1} & \num{1.399e-1} \\ \hline
  \end{tabular}
\end{adjustbox}
\end{center}
\end{table}

\begin{figure}[H]
  \centering
  \includegraphics[width=1\textwidth]{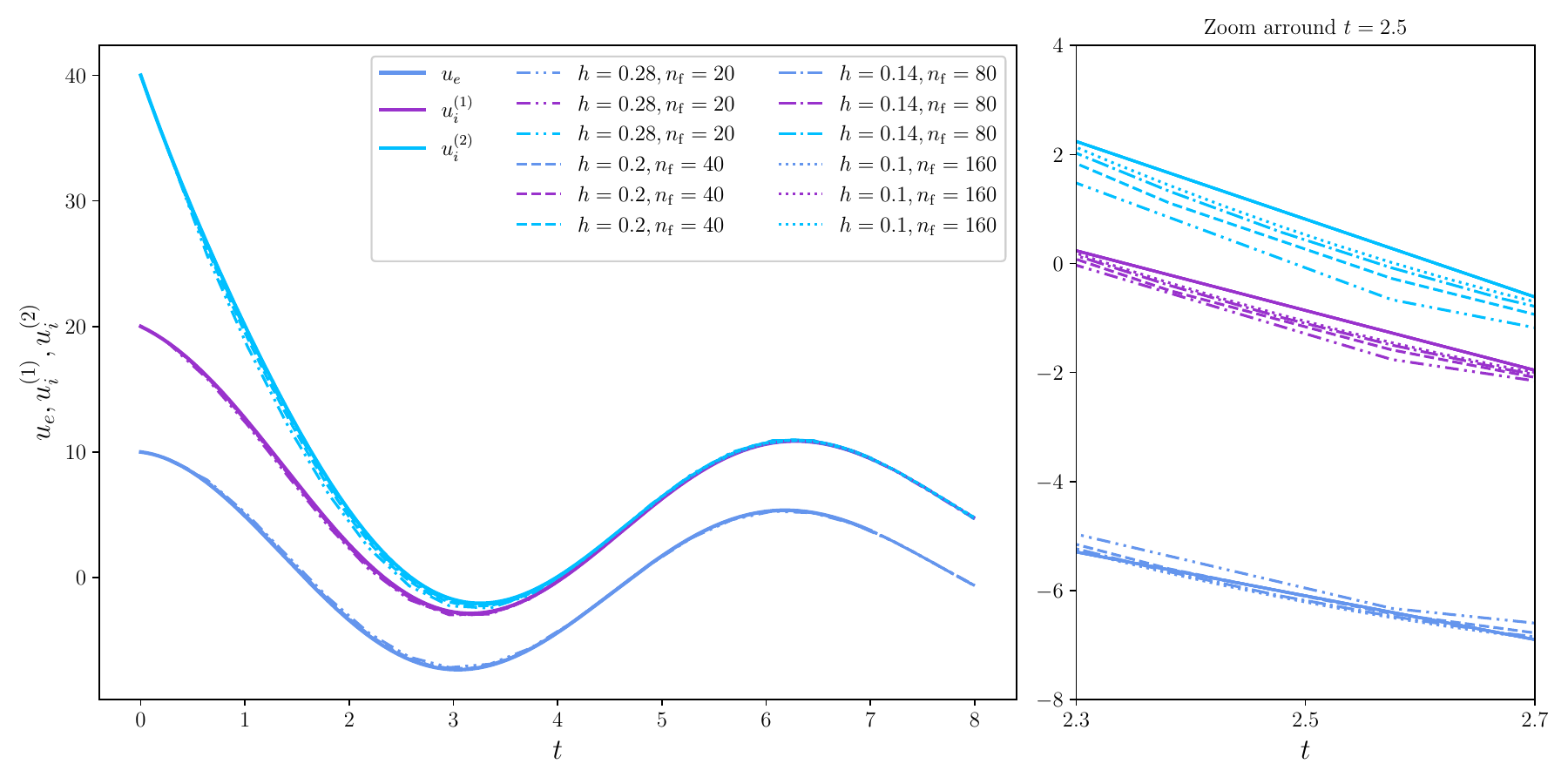}
  \caption{Numerical and exact solutions for $\ue$, $\ui[2]$, and
    $\ui[2]$ at $\rho=5$ in Experiment 2.}
    \label{fig:plot1_exp2}
\end{figure}

\begin{figure}[H]
  \centering
  \includegraphics[width=1\textwidth]{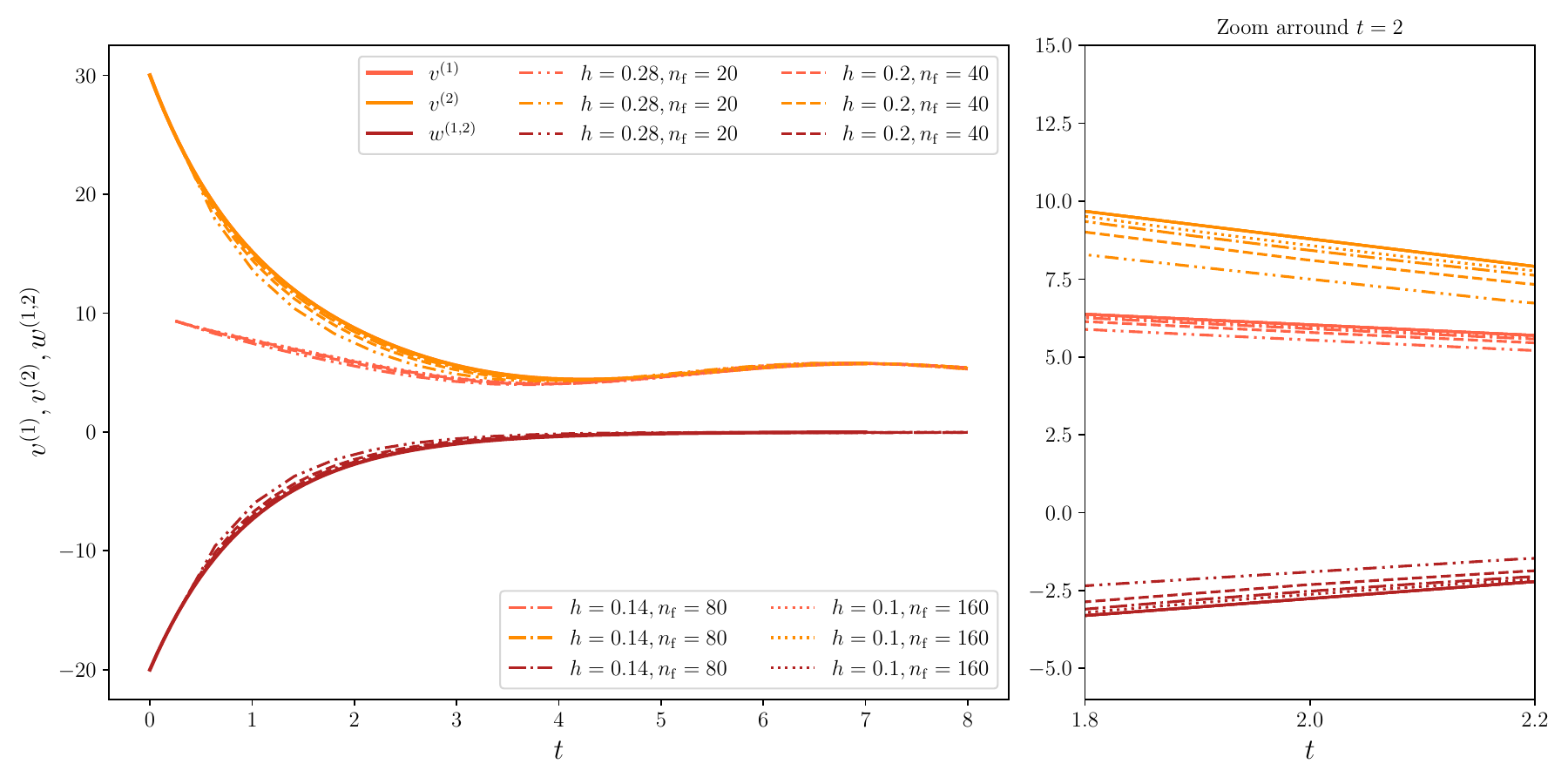}
  \caption{Numerical and exact solutions for $\vk[1]$ and $\vk[2]$ at
    $\rho=5$, and $\wk[1,2]$ at $\rho=4.5$ in Experiment 2.}
    \label{fig:plot2_exp2}
\end{figure}

In \cref{fig:plot1_exp2,fig:plot2_exp2}, we observe that 
the numerical solutions converge toward the exact solutions as 
the mesh size and time step decrease. Furthermore, \cref{tab:experiment-2} 
confirms the correct implementation of the numerical method. It is worth 
noting that the domain for the $L^2$-norm of the error in the variables $\ui[1]$, 
$\ui[2]$, and $\wk[1,2]$ is truncated to exclude points near the boundary. 
This is because Dirichlet boundary conditions are imposed, meaning 
that convergence near these points is not guaranteed due to the lack 
of control over the Neumann behavior.


\subsubsection{Experiment 3: 4-cell manufactured solution} 
\label{subsubsec:Experiment-3}
For the 4-cell manufactured solution, we consider the passive cell
model with parameters $\Rmk[1]=1$ and $\vrest=0$. Additionally, we set
the intracellular and extracellular conductivities to $\sigmai = 1$
and $\sigmae = 4$. For $\Te=1$, $\Ti[k]=k e^{-t}+\Te$ with cell
dimensions $\alpha_g=\beta_g=\gamma_g=0.5$ and $\alpha_c= \beta_c=\gamma_c=1$, 
the exact solution for the 4-cell system is given by

\begin{align*}
    \ue &= \cos(4 \pi x) \cos(4 \pi y) \cos(4 \pi z),\\
    \ui[k]&= \big(1 + k e^{-t}\big)\cos(4 \pi x) \cos(4 \pi y) \cos(4 \pi z),\\
    \vk[k]&= k e^{-t}\cos(4 \pi x) \cos(4 \pi y) \cos(4 \pi z),\\
    \wk[k,\ell]&= (k-\ell)e^{-t}\cos(4 \pi x) \cos(4 \pi y) \cos(4 \pi z),
\end{align*}
with forcing terms
\begin{align*}
    \fe &= 192 \pi^2 \cos(4 \pi x) \cos(4 \pi y) \cos(4 \pi z),\\
    \fik[k] &= 48\big(1 + ke^{-t}\big)\pi^2\cos(4 \pi x) \cos(4 \pi y) \cos(4 \pi z),\\
    \gk[k]&=\gk[k, \ell] = 0,
\end{align*}
for $ k, l \in \left\{1, 2, 3, 4 \right\}$. Enclosing the cells in a cube 
with dimensions $4.75\times4.75\times1.75$ to represent the extracellular 
space, we set the extracellular boundary condition to
\begin{align*}
    \uapp = 0.
\end{align*}

\setlength{\tabcolsep}{10pt}
\renewcommand{\arraystretch}{1.2}
\begin{table}[H]
  \caption{Errors for intracellular and extracelullar potentials in
    Experiment 3.  The simulation time spans are from $\tn[0]=0$ to
    $\tend=1$, $\nf$ denotes the number of time steps, and $\cl$ the
    characteristic length of the mesh.
  }\label{tab:experiment-3-intra-extra-spaces}
\begin{center}
\begin{adjustbox}{width=1 \textwidth, center}
  \begin{tabular}{c c c c c c c } \hline
   \rule{0pt}{18pt}
   $\cl$ & $\nf$ & 
   $\left\| \uet - \ue \right\|_{2}$ & 
   $\left\| \uit[1] - \ui[1] \right\|_{2}$ & 
   $\left\| \uit[2] - \ui[2] \right\|_{2}$ &
   $\left\| \uit[3] - \ui[3] \right\|_{2}$ &
   $\left\| \uit[4] - \ui[4] \right\|_{2}$ \\[7pt] \hline
    \num{  1.4e-1}  &  10  & \num{4.403e-1} & \num{1.342e-1} & \num{3.768e-1} & \num{1.601e-1} & \num{4.337e-1} \\
    \num{  1.0e-1}  &  20  & \num{2.902e-1} & \num{8.653e-2} & \num{2.400e-1} & \num{1.053e-1} & \num{2.869e-1} \\
    \num{  7.0e-2}  &  40  & \num{1.666e-1} & \num{4.197e-2} & \num{1.158e-1} & \num{5.348e-2} & \num{1.511e-1} \\
    \num{  4.9e-2}  &  80  & \num{9.349e-2} & \num{2.368e-2} & \num{6.612e-2} & \num{3.042e-2} & \num{8.897e-2} \\
    \num{  3.5e-2}  &  160 & \num{4.982e-2} & \num{1.295e-2} & \num{3.385e-2} & \num{1.660e-2} & \num{2.021e-2} \\ \hline
  \end{tabular}
\end{adjustbox}
\end{center}
\end{table}

\setlength{\tabcolsep}{10pt}
\renewcommand{\arraystretch}{1.2}
\begin{table}[H]
  \caption{Errors for transmembrane potentials in Experiment 3.  The
    simulation time spans are from $\tn[0]=0$ to $\tend=1$, $\nf$
    denotes the number of time steps, and $\cl$ the characteristic
    length of the mesh.  }\label{tab:experiment-3-membranes}
\begin{center}
\begin{adjustbox}{width=1 \textwidth, center}
  \begin{tabular}{c c c c c c } \hline
   \rule{0pt}{15pt}
   $\cl$ & $\nf$ & 
   $\left\| \vkt[1] - \vk[1] \right\|_{2}$ &
   $\left\| \vkt[2] - \vk[2] \right\|_{2}$ &
   $\left\| \vkt[3] - \vk[3] \right\|_{2}$ &
   $\left\| \vkt[4] - \vk[4] \right\|_{2}$ \\[5pt] \hline
    \num{  1.4e-1}  &  10  & \num{1.973e-1} & \num{5.414e-1} & \num{2.367e-1} & \num{6.351e-1} \\
    \num{  1.0e-1}  &  20  & \num{1.353e-1} & \num{3.711e-1} & \num{1.583e-1} & \num{4.414e-1} \\
    \num{  7.0e-2}  &  40  & \num{6.336e-2} & \num{1.867e-1} & \num{7.525e-2} & \num{2.408e-1} \\
    \num{  4.9e-2}  &  80  & \num{3.623e-2} & \num{1.104e-1} & \num{4.291e-2} & \num{1.400e-1} \\
    \num{  3.5e-2}  &  160 & \num{2.021e-2} & \num{6.181e-2} & \num{2.338e-2} & \num{7.824e-2}  \\ \hline
      \rule{0pt}{15pt}
    $\cl$ & $\nf$ & 
    $\left\| \wkt[1,2] - \wk[1,2] \right\|_{2}$ &
    $\left\| \wkt[1,3] - \wk[1,3] \right\|_{2}$ &
    $\left\| \wkt[2,4] - \wk[2,4] \right\|_{2}$ &
    $\left\| \wkt[3,4] - \wk[3,4] \right\|_{2}$ \\[5pt] \hline
    \num{  1.4e-1}  &  10  & \num{2.821e-2} & \num{4.395e-2} & \num{5.263e-2} & \num{7.484e-2} \\
    \num{  1.0e-1}  &  20  & \num{3.944e-2} & \num{3.376e-2} & \num{4.871e-2} & \num{5.011e-2} \\
    \num{  7.0e-2}  &  40  & \num{1.904e-2} & \num{2.596e-2} & \num{2.700e-2} & \num{2.941e-2} \\
    \num{  4.9e-2}  &  80  & \num{1.111e-2} & \num{1.375e-2} & \num{1.546e-2} & \num{1.711e-2} \\
    \num{  3.5e-2}  &  160 & \num{5.800e-3} & \num{6.865e-3} & \num{8.464e-3} & \num{9.078e-3} \\ \hline
  \end{tabular}
\end{adjustbox}
\end{center}
\end{table}

The results presented in
\Cref{tab:experiment-3-intra-extra-spaces,,tab:experiment-3-membranes}
validate the correct implementation of the numerical method and the
accuracy of the 4-cell manufactured solution. The $L^2$-norm of the
errors for the variables $\ue$, $\vk[k]$, $\ui[k]$, and $\wk[k,\ell]$
with $k, l \in \left\{1, 2, 3, 4\right\}$ consistently decrease as the mesh size
and time step are refined.

\section{Conclusions}
\label{sec:Conclusions}

The EMI model provides a powerful framework for modeling cardiac
tissue.  However, its high computational cost limits its application.
This limitation drives the development of new computational methods
and hence highlights the utility of analytical solutions for
establishing the accuracy and efficiency of software.  In this work,
we derived a family of analytical solutions for the EMI model for both
a single circular cell (\cref{eq:sol-single-cell}) and two coupled
semi-spherical cells (\cref{eq:sol-two-coupled-cells}), assuming
radial symmetry and equal conductivities in the intracellular and
extracellular spaces.  This last assumption has relevance when
modeling liposomes.  Although these solutions diverge at the origin,
we demonstrate that they can be accurately simulated by excluding the
origin from the domain and imposing Dirichlet boundary conditions on
the resulting boundary for the corresponding variables.  This approach
makes the solutions suitable for the convergence analysis of numerical
methods. We simulated some simple examples to demonstrate the
convergence of numerical solutions to the analytical ones; see,
\cref{subsubsec:Experiment-1,subsubsec:Experiment-2}.  Additionally,
we present a manufactured solution for $N$ three-dimensional cells
(\cref{eq:sol-manufactured}) that was used to simulate 4 cells and
validated the accuracy of the numerical implementation; see,
\cref{subsubsec:Experiment-3}.  Altogether, these solution form a
benchmark for assessing the accuracy of numerical methods for solving
the EMI model.  Furthermore, the Python implementation of the
numerical method provides a solid foundation for future research on
the convergence analysis of high-order operator splitting methods.

\section*{Data availability}  
The dataset supporting this study is available at
\url{https://doi.org/10.5281/zenodo.14948002}.  The source code to
reproduce the data can be accessed in the following repository:
\url{https://github.com/uofs-simlab/EMI-Analytical-Solutions}.



\bibliographystyle{elsarticle-num} 
\bibliography{references}

\end{document}